\documentclass[draft]{amsart}
\usepackage{mathrsfs}
\usepackage{amssymb}
\usepackage{amsmath}
\usepackage{amsfonts}
\usepackage{amsfonts,enumerate}
\usepackage{pifont}
\usepackage{enumerate}
\usepackage{graphics}
\usepackage{verbatim}
\usepackage{amsthm}

\numberwithin{equation}{section}
\newtheorem{theorem}{Theorem}[section]
\newtheorem{corollary}{Corollary}[section]

\newtheorem{lemma}{Lemma}[section]
\newtheorem{proposition}{Proposition}[section]
\newtheorem{remark}{Remark}[section]

\newcommand{\A}{{\mathcal A}}

\newcommand{\g}{\gamma}

\newcommand{\8}{\infty}
\newcommand{\el}{\ell}

\newcommand{\be}{\begin{eqnarray*}}
\newcommand{\ee}{\end{eqnarray*}}
\newcommand{\beq}{\begin{equation}}
\newcommand{\eeq}{\end{equation}}
\newcommand{\beqn}{\begin{equation*}}
\newcommand{\eeqn}{\end{equation*}}
\newcommand{\bs}{\begin{split}}
\newcommand{\es}{\end{split}}

\numberwithin{equation}{section}

\begin{document}

\title{Maximal and area integral characterizations of Bergman spaces in the unit ball of $\mathbb{C}^n$}


\author{Zeqian Chen}

\address{Wuhan Institute of Physics and Mathematics, Chinese
Academy of Sciences, 30 West District, Xiao-Hong-Shan, Wuhan 430071, China}


\author{Wei Ouyang}

\address{Wuhan Institute of Physics and Mathematics, Chinese
Academy of Sciences, 30 West District, Xiao-Hong-Shan, Wuhan 430071, China
and Graduate University of Chinese Academy of Sciences, Beijing 100049, China}


\date{}
\maketitle

\markboth{Z. Chen and W. Ouyang}%
{Bergman spaces}

\begin{abstract}
In this paper, we present maximal and area integral characterizations of Bergman spaces in the unit ball of $\mathbb{C}^n.$ The characterizations are in terms of maximal functions and area integral functions on Bergman balls involving the radial derivative, the complex gradient, and the invariant gradient.
As an application, we obtain new maximal and area integral characterizations of Besov spaces. Moreover, we give an atomic decomposition of real-variable type with respect to Carleson tubes for Bergman spaces.
\end{abstract}


\section{Introduction and main results}\label{intro}

Let $\mathbb{C}$ denote the set of complex numbers. Throughout the paper we fix a positive integer $n,$ and let
\be
\mathbb{C}^n=\mathbb{C}\times\cdots\times\mathbb{C}
\ee
denote the Euclidean space of complex dimension $n.$ Addition, scalar multiplication, and conjugation are defined on $\mathbb{C}^n$ componentwise. For $z=(z_1,\cdots,z_n)$ and $w=(w_1,\cdots,w_n)$ in $\mathbb{C}^n,$ we write
\be
\langle z,w\rangle=z_1\overline{w}_1+\cdots+z_n\overline{w}_n,
\ee
where $\overline{w}_k$ is the complex conjugate of $w_k.$ We also write
\be
|z|=\sqrt{|z_1|^2+\cdots+|z_n|^2}.
\ee
The open unit ball in $\mathbb{C}^n$ is the set
\be
{\mathbb{B}_n} =\{\,z\in {\mathbb{C}^n}: |z|<1\}.
\ee
The boundary of $\mathbb{B}_n$ will be denoted by $\mathbb{S}_n$ and is called the unit sphere in $\mathbb{C}^n,$ i.e.,
\be
{\mathbb{S}_n}=\{\,z\in {\mathbb{C}^n}: |z|=1\}.
\ee
Also, we denote by $\overline{\mathbb{B}}_n$ the closed unit ball, i.e.,
\be
\overline{\mathbb{B}}_n= \{z \in \mathbb{C}^n:\; |z| \le 1 \} = \mathbb{B}_n \cup \mathbb{S}_n.
\ee
The automorphism group of $\mathbb{B}^n,$ denoted by $\mathrm{Aut} (\mathbb{B}^n),$ consists of all bi-holomorphic mappings of
$\mathbb{B}^n.$ Traditionally, bi-holomorphic mappings are also called automorphisms.

For $\alpha \in \mathbb{R},$ the weighted Lebesgue measure $dv_{\alpha}$ on $\mathbb{B}_n$ is defined by
\be
dv_{\alpha}(z)=c_{\alpha}(1-|z|^2)^{\alpha}dv(z)
\ee
where $c_{\alpha}=1$ for $\alpha\le-1$ and $c_{\alpha}=\Gamma(n+\alpha+1)/ [n!\Gamma(\alpha+1)]$ if $\alpha>-1$, which is a normalizing constant so that $dv_{\alpha}$
is a probability measure on $\mathbb{B}_n.$ In the case of $\alpha=-(n+1)$ we denote the resulting measure by
\be
d\tau(z)=\frac{dv}{(1-|z|^2)^{n+1}},
\ee
and call it the invariant measure on $\mathbb{B}^n,$ since $d\tau=d\tau\circ\varphi$ for any automorphism $\varphi$ of $\mathbb{B}^n.$

For $\alpha > -1$ and $p>0,$ the (weighted) Bergman space $\mathcal{A}^p_{\alpha}$ consists of holomorphic functions $f$ in $\mathbb{B}_n$ with
$$\|f\|_{p,\,\alpha}=\left ( \int_{\mathbb{B}_n}|f(z)|^pdv_{\alpha}(z) \right )^{1/p}<\infty,$$ where the weighted Lebesgue measure $dv_{\alpha}$ on $\mathbb{B}_n$ is defined by
$$dv_{\alpha}(z)=c_{\alpha}(1-|z|^2)^{\alpha}dv(z)$$
and $c_{\alpha}=\Gamma(n+\alpha+1)/ [n!\Gamma(\alpha+1)]$ is a normalizing constant so that $dv_{\alpha}$ is a probability measure on $\mathbb{B}_n.$ Thus,
\be
\mathcal{A}^p_{\alpha} = \mathcal{H} (\mathbb{B}_n) \cap L^p (\mathbb{B}_n, d v_{\alpha}),
\ee
where $\mathcal{H} (\mathbb{B}_n)$ is the space of all holomorphic functions in $\mathbb{B}_n.$ When $\alpha =0$ we simply write $\A^p$ for $\A^p_0.$ These are the usual Bergman spaces. Note that for $1 \le p < \8,$ $\mathcal{A}^p_{\alpha}$ is a Banach space under the norm $\|\ \|_{p,\,\alpha}.$ If $0 < p <1,$ the space $\mathcal{A}^p_{\alpha}$ is a quasi-Banach space with $p$-norm $\| f \|^p_{p, \alpha}.$

Recall that $D(z,\gamma)$ denotes the Bergman metric ball at $z$
\be
D(z, \gamma) = \{w \in \mathbb{B}_n\;: \; \beta (z, w) < \g \}
\ee
with $\gamma >0,$ where $\beta$ is the Bergman metric on $\mathbb{B}_n.$ It is known that
\be
\beta (z, w) = \frac{1}{2} \log \frac{1 + | \varphi_z (w)|}{1 - | \varphi_z (w)|},\quad z, w \in \mathbb{B}_n,
\ee
whereafter $\varphi _z$ is the bijective holomorphic mapping in $\mathbb{B}_n,$ which satisfies $\varphi _z (0)=z$, $\varphi _z(z)=0$ and
$\varphi _z\circ\varphi _z = id.$

As is well known, maximal functions play a crucial role in the real-variable theory of Hardy spaces (cf. \cite{Stein1993}). In this paper, we first establish a maximal-function characterization for the Bergman spaces. To this end, we define for each $\g > 0$ and $f \in \mathcal{H} (\mathbb{B}_n):$
\beq\label{eq:MaximalFunct}
(\mathrm{M}_{\g} f) (z) = \sup_{w \in D(z, \g)} | f (w)|,\; \forall z \in \mathbb{B}_n.
\eeq
Then we have the following result.

\begin{theorem}\label{th:MaxialCharat}
Suppose $\gamma >0$ and $\alpha > -1.$ Let $0< p < \8.$ Then for any $f \in \mathcal{H} (\mathbb{B}_n),$ $f \in \mathcal{A}^p_{\alpha}$ if and only if $\mathrm{M}_{\g} f \in L^p (\mathbb{B}_n, d v_{\alpha}).$ Moreover,
\beq\label{eq:MaximalFunctNorm}
\| f \|_{p,\alpha} \approx \| \mathrm{M}_{\g} f \|_{p, \alpha},
\eeq
where ``$\approx$" depends only on $\gamma, \alpha, p,$ and $n.$
\end{theorem}

The norm appearing on the right-hand side of \eqref{eq:MaximalFunctNorm} can be viewed an analogue of the so-called nontangential maximal function in Hardy spaces. The proof of Theorem \ref{th:MaxialCharat} is fairly elementary (see \S \ref{PrCmethod}), using some basic facts and estimates on the Bergman balls.

In order to state the area integral characterizations of the Bergman spaces, we require some more notation. For any $f \in \mathcal{H} ( \mathbb{B}_n)$ and $z = (z_1, \ldots, z_n) \in \mathbb{B}_n$ we define
\be
\mathcal{R} f (z) = \sum^n_{k=1} z_k \frac{\partial f (z)}{\partial z_k}
\ee
and call it the radial derivative of $f$ at $z.$ The complex and invariant gradients of $f$ at $z$ are respectively defined as
\be
\nabla f(z) = \Big ( \frac{\partial f (z) }{\partial z_1},\ldots, \frac{\partial f (z) }{\partial z_n} \Big )\; \text{and}\; \widetilde{\nabla} f (z) = \nabla(f \circ \varphi_z)(0).
\ee

Now, for fixed $1< q < \8$ and $\gamma >0,$ we define for each $f \in \mathcal{H} (\mathbb{B}_n)$ and $z \in \mathbb{B}_n:$
\begin{enumerate}[{\rm (1)}]

\item The radial area integral function
\be
A_{\mathcal{R}}^{\gamma, q} (f) (z) = \left ( \int_{D(z,\gamma)}  | (1-|w|^2) \mathcal{R} f (w) |^q d \tau (w) \right )^{\frac{1}{q}}.
\ee

\item The complex gradient area integral function
\be
A_{\nabla}^{\gamma, q} ( f) (z) = \left ( \int_{D(z, \gamma)} | (1-|w|^2) \nabla f (w) |^q d \tau (w) \right )^{\frac{1}{q}}.
\ee

\item The invariant gradient area integral function
\be
A_{\tilde{\nabla}}^{\gamma, q} (f) (z) = \left ( \int_{D(z, \gamma)} | \tilde{\nabla} f (w) |^q d \tau (w) \right )^{\frac{1}{q}}.
\ee

\end{enumerate}

We state the second main result of this paper as follows.

\begin{theorem}\label{th:AreaCharat}
Suppose $1< q < \8, \gamma > 0,$ and $\alpha > -1.$ Let $0 < p < \8.$ Then, for any $f \in \mathcal{H} (\mathbb{B}_n)$ the following conditions are equivalent:
\begin{enumerate}[{\rm (a)}]

\item $f \in \mathcal{A}^p_{\alpha}.$

\item $A_{\mathcal{R}}^{\gamma, q} (f)$ is in $L^p (\mathbb{B}_n, d v_{\alpha}).$

\item $A_{\nabla}^{\gamma, q} ( f)$ is in $L^p (\mathbb{B}_n, d v_{\alpha}).$

\item $A_{\tilde{\nabla}}^{\gamma, q} (f)$ is in $L^p (\mathbb{B}_n, d v_{\alpha}).$

\end{enumerate}
Moreover, the quantities
\be
\| A_{\mathcal{R}}^{\gamma, q} (f) \|_{p, \alpha},\; \| A_{\nabla}^{\gamma, q} ( f) \|_{p, \alpha}, \; \| A_{\tilde{\nabla}}^{\gamma, q} (f) \|_{p, \alpha},
\ee
are all comparable to $\| f - f(0) \|_{p, \alpha},$ where the comparable constants depend only on $q, \gamma, \alpha, p,$ and $n.$
\end{theorem}

For $0 < p < \8$ and $-\8 < \alpha < \8,$ we fix a nonnegative integer $k$ with $pk + \alpha > -1$ and define the generalized Bergman space $\mathcal{A}^p_{\alpha}$ as introduced in \cite{ZZ2008} to be the space of all $f \in \mathcal{H} (\mathbb{B}_n)$ such that $(1-|z|^2)^k \mathcal{R}^k f \in L^p (\mathbb{B}_n, d v_{\alpha}).$ One then easily observes that $\mathcal{A}^p_{\alpha}$ is independent of the choice of $k$ and consistent with the traditional definition when $\alpha > -1.$ Let $N$ be the smallest nonnegative integer
such that $pN + \alpha > -1.$ Put
\beq\label{eq:BergmanSpaceNorm}
\| f \|_{p, \alpha} = | f(0)| + \left ( \int_{\mathbb{B}_n} (1-|z|^2)^{pN} | \mathcal{R}^N f (z) |^p d v_{\alpha} (z) \right )^{\frac{1}{p}},\quad f \in \mathcal{A}^p_{\alpha}\;.
\eeq
Equipped with \eqref{eq:BergmanSpaceNorm}, $\mathcal{A}^p_{\alpha}$ becomes a Banach space when $p \ge 1$ and a quasi-Banach space for $0 < p < 1.$

It is known that the family of the generalized Bergman spaces $\mathcal{A}^p_{\alpha}$ covers most of the spaces of holomorphic functions in the unit ball of $\mathbb{C}^n,$ such as the classical diagonal Besov space $B^s_p$ and the Sobolev space $W^p_{k,\beta},$ which has been extensively studied before in the literature under different names (see e.g. \cite{ZZ2008} for an overview).

There are various characterizations for $B^s_p$ or $W^p_{k,\beta}$ involving complex-variable quantities in terms of radical derivatives, complex and invariant gradients, and fractional differential operators (for a review and details see \cite{ZZ2008} and references therein). However, as an application of Theorems \ref{th:MaxialCharat} and \ref{th:AreaCharat}, we obtain new maximal and area integral characterizations of the Besov spaces as follows, which can be considered as a unified characterization for such spaces involving real-variable quantities.

\begin{corollary}\label{cor:MaxialCharat}
Suppose $\gamma >0$ and $\alpha \in \mathbb{R}.$ Let $0< p < \8$ and $k$ be a positive integer such that $pk + \alpha > -1.$ Then for any $f \in \mathcal{H} (\mathbb{B}_n),$ $f \in \mathcal{A}^p_{\alpha}$ if and only if $\mathrm{M}^{(k)}_{\g} ( f ) \in L^p (\mathbb{B}_n, d v_{\alpha}),$
where
\beq\label{eq:kMaximalFunct}
\mathrm{M}^{(k)}_{\g} ( f ) (z) = \sup_{w \in D (z, \gamma)} | (1-|w|^2)^k \mathcal{R}^k f (w) |,\quad z \in \mathbb{B}_n.
\eeq
Moreover,
\beq\label{eq:kMaximalFunctNorm}
\| f - f(0) \|_{p,\alpha} \approx \| \mathrm{M}_{\g} ( \mathcal{R}^k f ) \|_{p, \alpha},
\eeq
where ``$\approx$" depends only on $\gamma, \alpha, p, k,$ and $n.$
\end{corollary}

\begin{corollary}\label{cor:AreaCharat}
Suppose $1< q < \8, \gamma >0$ and $\alpha \in \mathbb{R}.$ Let $0 < p < \8$ and $k$ be a nonnegative integer such that $pk + \alpha > -1.$ Then for any $f \in \mathcal{H} (\mathbb{B}_n),$ $f \in \mathcal{A}^p_{\alpha}$ if and only if $A_{\mathcal{R}^{k+1}}^{\gamma, q} ( f)$ is in $L^p (\mathbb{B}_n, d v_{\alpha}),$
where
\beq\label{eq:kareaFunct}
A_{\mathcal{R}^{k+1}}^{\gamma, q} ( f) (z) = \left ( \int_{D(z,\gamma)} \big | (1-|w|^2)^{k+1} \mathcal{R}^{k+1} f (w) \big |^q d \tau (w) \right )^{\frac{1}{q}}.
\eeq
Moreover,
\beq\label{eq:kAreaFunctNorm}
\| f - f(0) \|_{p,\alpha} \approx \| A_{\mathcal{R}^{k+1}}^{\gamma, q} ( f) \|_{p, \alpha},
\eeq
where ``$\approx$" depends only on $q, \gamma, \alpha, p, k,$ and $n.$
\end{corollary}

To prove Corollaries \ref{cor:MaxialCharat} and \ref{cor:AreaCharat}, one merely notices that $f \in \mathcal{A}^p_{\alpha}$ if and only if $\mathcal{R}^k f \in L^p (\mathbb{B}_n, d v_{\alpha + pk})$ and applies Theorems \ref{th:MaxialCharat} and \ref{th:AreaCharat} respectively to $\mathcal{R}^k f$ with the help of Lemma \ref{le:KerEstimation} below.

The paper is organized as follows. In Sect. \ref{PrCmethod} we will prove Theorems \ref{th:MaxialCharat} and \ref{th:AreaCharat}. An atomic decomposition of real-variable type with respect to Carleson tubes for Bergman spaces will be presented in Sect. \ref{Atomdecomp} via duality method. Finally, in Sect. \ref{AreaCharactp=1}, we will prove Theorem \ref{th:AreaCharat} through using the real-variable atomic decomposition of Bergman spaces established in the preceding section.

In what follows, $C$ always denotes a constant depending (possibly) on $n, q, p, \g$ or $\alpha$ but not on $f,$ which may be
different in different places. For two nonnegative (possibly infinite) quantities $X$ and $Y,$ by $X \lesssim Y$ we mean that there
exists a constant $C>0$ such that $ X \leq C Y$ and by $X \thickapprox Y$ that $X \lesssim Y$ and $Y \lesssim X.$ Any notation and terminology not otherwise explained, are as used in \cite{Zhu2005} for spaces of holomorphic functions in the unit ball of $\mathbb{C}^n.$

\section{Proofs of Theorems \ref{th:MaxialCharat} and \ref{th:AreaCharat}}\label{PrCmethod}

For the sake of convenience, we collect some elementary facts on the Bergman metric and holomorphic functions in the unit ball of $\mathbb{C}^n$ as follows.

\begin{lemma}\label{le:KerEstimation} {\rm (cf. \cite[Lemma 2.20]{Zhu2005})}
For each $\gamma>0,$
\be
1-|a|^2 \approx 1-|z|^2 \approx |1-\langle a,z\rangle|
\ee
for all $a$ and $z$ in $\mathbb{B}_n$ with $\beta(a,z)<{\gamma}.$
\end{lemma}

\begin{lemma}\label{le:PointEstimation} {\rm (cf. \cite[Lemma 2.24]{Zhu2005})}
Suppose $\g >0, p>0,$ and $\alpha > -1.$ Then there exists a constant $C>0$ such that for any $f \in \mathcal{H} (\mathbb{B}_n),$
\be
|f(z) |^p \le \frac{C}{(1-|z|^2)^{n+1+\alpha}} \int_{D(z,\g)} | f(w)|^p d v_{\alpha} (w), \quad \forall z \in \mathbb{B}_n.
\ee
\end{lemma}

\begin{lemma}\label{le:KerEstimation2} {\rm (cf. \cite[Lemma 2.27]{Zhu2005})}
For each $\gamma>0,$
\be
|1-\langle z,u\rangle| \approx |1-\langle z,v\rangle|
\ee
for all $z$ in $\bar{\mathbb{B}}_n$ and $u,v$ in $\mathbb{B}_n$ with $\beta(u,v)<\gamma.$
\end{lemma}

\subsection{Proof of Theorem \ref{th:MaxialCharat}}

We need the following result (cf. \cite[Lemma 5]{BBCZ1990}).

\begin{lemma}\label{le:CoverBall}
For fixed $\gamma >0,$ there exist a positive integer $N$ and a sequence $\{a_k\}$ in $\mathbb{B}_n$ such that
\begin{enumerate}[{\rm (1)}]

\item $\mathbb{B}_n = \cup_k D(a_k, \gamma),$ and

\item each $z \in \mathbb{B}_n$ belongs to at most $N$ of the sets $D (a_k, 3 \gamma).$

\end{enumerate}
\end{lemma}

\

{\it Proof of Theorem \ref{th:MaxialCharat}}.\;
Let $p >0.$ By Lemmas \ref{le:CoverBall}, \ref{le:PointEstimation}, and \ref{le:KerEstimation},
we have
\be\begin{split}
\int_{\mathbb{B}_n} | \mathrm{M}_{\g}(f)(z)|^{p} & d v_{\alpha}(z)
\leq \sum_k \int_{D(a_{k},\gamma)} |\mathrm{M}_{\g} (f)(z)|^pdv_{\alpha} (z)\\
=&\sum_k \int_{D(a_{k},\gamma)}\sup_{w\in D(z,\gamma)}|f(w)|^pdv_{\alpha}(z)\\
\lesssim & \sum_k \int_{D(a_{k},\gamma)} \sup_{w\in D(z,\gamma)} \frac{1}{(1-|w|^2)^{n+1+\alpha}} \int_{D(w,\gamma)}|f(u)|^pdv_{\alpha}(u) d v_{\alpha}(z)\\
\lesssim & \sum_k \int_{D(a_{k},\gamma)} \left ( \frac{1}{(1-|a_{k}|^2)^{n+1+\alpha}}\int_{D(a_{k}, 3 \gamma)}|f(u)|^pdv_{\alpha}(u) \right ) d v_{\alpha}(z)\\
\lesssim & \sum_k \int_{D(a_{k}, 3 \gamma)}|f(u)|^pdv_{\alpha}(u)
\lesssim N \int_{\mathbb{B}_n}|f(u)|^pdv_{\alpha}(u)
\end{split}\ee
where $N$ is the constant in Lemma \ref{le:CoverBall} depending only on $\gamma$ and $n.$
\hfill $\Box$

\subsection{Proof of Theorem \ref{th:AreaCharat}}

Recall that $\mathcal{B} (\mathbb{B}_n)$ is defined as the space of all $f \in \mathcal{H}(\mathbb{B}_n)$ so that
\be
\| f \|_{\mathcal{B}} = \sup_{z \in \mathbb{B}_n} | \tilde{\nabla} f (z) | < \8.
\ee
$\mathcal{B} (\mathbb{B}_n)$ with the norm $\| f \| =  | f(0) | + \| f \|_{\mathcal{B}}$ is a Banach space and called the Bloch space. Then, the following interpolation result holds.

\begin{lemma}\label{le:InterpBloch}{\rm (cf. \cite[Theorem 3.25]{Zhu2005})}
Let $1 < p < \8.$ Suppose $\alpha>-1$ and
\be
\frac{1}{p}=\frac{1-\theta}{p'}
\ee
for $0 < \theta < 1$ and $1 \le p' < \8.$ Then
\be
\mathcal{A}^p_{\alpha}(\mathbb{B}_n) = \left [ \mathcal{A}^{p'}_{\alpha}(\mathbb{B}_n), \mathcal{B} (\mathbb{B}_n) \right ]_{\theta}
\ee
with equivalent norms.
\end{lemma}

Moreover, to prove Theorem \ref{th:AreaCharat} for the case $0< p \le 1,$ we will use atom decomposition for Bergman spaces due to Coifman and Rochberg \cite{CR1980} (see also \cite[Theorem 2.30]{Zhu2005}) as follows.

\begin{proposition}\label{prop:AtomDecpCR}
Suppose $p>0, \alpha>-1,$ and $b> n \max \{1, 1/p \} + (\alpha+1)/p.$ Then there exists a sequence $\{a_k\}$ in $\mathbb{B}_n$ such that
$\mathcal{A}^p_{\alpha}$ consists exactly of functions of the form
\be
f(z)=\sum^{\infty}_{k=1}c_{k}\frac{(1-|a_k|^2)^{(pb-n-1-\alpha)/p}}{(1-\langle z,a_{k}\rangle)^b},\ \ \  z\in\mathbb{B}_n,
\ee
where $\{c_k\}$ belongs to the sequence space $\el^p$ and the series converges in the norm topology of $\mathcal{A}^p_{\alpha}.$ Moreover,
\be
\int_{\mathbb{B}_n}|f(z)|^pdv_{\alpha}(z) \approx  \inf \Big \{ \sum_{k}|c_{k}|^p \Big \},
\ee
where the infimum runs over all the above decompositions.
\end{proposition}

Also, we need a characterization of Carleson type measures for Bergman spaces as follows, which can be found in \cite[Theorem 45]{ZZ2008}.

\begin{proposition}\label{prop:CarlesonMeasure}
Suppose $n+1+\alpha>0$ and $\mu$ is a positive Borel measure on $\mathbb{B}_n.$ Then, there exists a constant $C>0$ such that
\be
\mu(Q_{r}(\zeta))\le Cr^{2(n+1+\alpha)},\quad \forall \zeta\in{\mathbb{S}_n} \;\text{and}\;r>0,
\ee
if and only if for each $s>0$ there exists a constant $C>0$ such that
\be
\int_{\mathbb{B}_n}\frac{(1-|z|^2)^s}{|1-\langle z, w\rangle|^{n+1+\alpha+s}}d\mu(w) \le C
\ee
for all $z\in{\mathbb{B}_n}.$
\end{proposition}

We are now ready to prove Theorem \ref{th:AreaCharat}. Note that for any $f \in \mathcal{H} (\mathbb{B}_n),$
\be
(1-|z|^2) | \mathcal{R} f (z) | \le (1-|z|^2) |\nabla f (z)| \le | \tilde{\nabla} f (z) |, \quad \forall z \in \mathbb{B}_n
\ee
(cf. \cite[Lemma 2.14]{Zhu2005}). We have that (d) implies (c), and (c) implies (b) in Theorem \ref{th:AreaCharat}. Then, it remains to prove that (b) implies (a), and (a) implies (d).

\

{\it Proof of $(\mathrm{b}) \Rightarrow (\mathrm{a})$}.\;  Since $\mathcal{R} f(z)$ is holomorphic, by Lemma \ref{le:PointEstimation} we
have
\be\begin{split}
| \mathcal{R} f(z)|^q &\leq \frac{C}{(1-|z|^2)^{n+1}} \int_{D(z,\gamma)}| \mathcal{R} f(w)|^q d v(w) \leq C_{\g} \int_{D(z,\gamma)} | \mathcal{R} f(w)|^q d \tau (w).
\end{split}\ee
Then,
\be\begin{split}
(1-|z|^2) | \mathcal{R} f (z)|\leq & C(1-|z|^2)\left( \int_{D(z,\gamma)} | \mathcal{R} f(w)|^q d \tau (w) \right)^{\frac{1}{q}}\\
\leq & C_{\g} \left ( \int_{D(z,\gamma)}| (1-|w|^2) \mathcal{R} f(w)|^q d \tau (w) \right )^{\frac{1}{q}}
= C_{\g} A^{\gamma, q}_{\mathcal{R}} (f) (z).
\end{split}\ee
Hence, for any $p >0,$ if $A^{\gamma, q}_{\mathcal{R}} (f) \in L^p (\mathbb{B}_n, dv_{\alpha})$ then $(1-|z|^2)| \mathcal{R} f(z)|$ is in $L^p (\mathbb{B}_n, d v_{\alpha}),$ which implies that $f \in \mathcal{A}^p_{\alpha}$ (cf. \cite[Theorem 2.16]{Zhu2005}).
\hfill$\Box$

\

The proof of $(\mathrm{a}) \Rightarrow (\mathrm{d})$ is divided into two steps. We first prove the case $0< p \le 1$ using the atomic decomposition, and then the remaining case via complex interpolation.

{\it Proof of $(\mathrm{a}) \Rightarrow (\mathrm{d})$ for $0 < p \le 1$}.\;
To this end, we write
\be
f_k(z)=\frac{(1-|a_k|^2)^{(pb-n-1-\alpha)/p}}{(1-\langle z,a_{k}\rangle)^b}.
\ee
An immediate computation yields that
\be
\nabla f_k(z)=\frac{ b \overline{a}_k(1-|a_k|^2)^{(pb-n-1-\alpha)/p}}{(1-\langle
z,a_{k}\rangle)^{b+1}}
\ee
and
\be
\mathcal{R} f_{k}(z)=\frac{b\langle z,a_k\rangle(1-|a_k|^2)^{(pb-n-1-\alpha)/p}}{(1-\langle
z,a_{k}\rangle)^{b+1}}
\ee
Then we have
\be\begin{split}
|\tilde{\nabla} f_k(z)|^2=&(1-|z|^2)(| \nabla f_k (z) |^2 - | \mathcal{R} f_k (z)|^2 )\\
=& b^2(1-|z|^2)(1-|a_k|^2)^{2(pb-n-1-\alpha)/p}\frac{|a_k|^2-|\langle z,a_{k} \rangle|^2}{|1-\langle
z, a_{k}\rangle|^{2(b+1)}}.
\end{split}\ee
By Lemmas \ref{le:KerEstimation} and \ref{le:KerEstimation2} one has
\be\begin{split}
A^{\gamma, q}_{\tilde{\nabla}}(f_{k})(z)
& = \left ( \int_{D(z, \gamma)}  | \tilde{\nabla} f_k (w) |^q d \tau (w) \right )^{\frac{1}{q}}\\
& \le b(1-|a_k|^2)^{(pb-n-1-\alpha)/p} \Big ( \int_{D(z, \gamma)} \frac{1}{|1-\langle w,a_k\rangle|^{q b}} d \tau (w) \Big )^{\frac{1}{q}}\\
& \le \frac{C_{\gamma} b(1-|a_k|^2)^{(pb-n-1-\alpha)/p}}{{|1-\langle z,a_k\rangle|^{b}}},
\end{split}\ee
where we have used the fact $v(D(z, \gamma))\approx (1-|z|^2)^{n+1}.$ Note that $v_{\alpha} (Q_r) \approx r^{2 (n+1 + \alpha)}$ (cf. \cite[Corollary 5.24]{Zhu2005}), by Proposition \ref{prop:CarlesonMeasure} we have
\be
\int_{\mathbb{B}_n}| A^{\gamma, q}_{\tilde{\nabla}}(f_{k}) (z)|^p d v_{\alpha} (z) \le C b^{p} \int_{\mathbb{B}_n}\frac{(1-|a_k|^2)^{(pb-n-1-\alpha)}}{{|1-\langle
z,a_k\rangle|^{pb}}}dv_{\alpha}(z)\le C_{p, \alpha}.
\ee
Hence, for $0 < p \le 1$ we have for $f = \sum^{\infty}_{k=1}c_{k} f_k$ with $\sum_k |c_k|^p < \8,$
\be\begin{split}
\int_{\mathbb{B}_n}| A^{\gamma, q}_{\tilde{\nabla}}(f) (z)|^p d v_{\alpha} \leq \sum_{k=1}^{\infty} |c_k|^p \int_{\mathbb{B}_n} | A^{\gamma, q}_{\tilde{\nabla}}(f_{k}) (z)|^pdv_{\alpha} \leq C_{p,\alpha} \sum_{k=1}^{\infty} | c_{k}|^{p}.
\end{split}\ee
This concludes that
\be
\int_{\mathbb{B}_n}| A^{\gamma, q}_{\tilde{\nabla}}(f) (z)|^p dv_{\alpha}\leq C_{p, \alpha} \inf \Big \{ \sum_{k=1}^{\infty} |c_{k}|^{p} \Big \} \leq C_{p, \alpha} \int_{\mathbb{B}_n} |f(z)|^p d v_{\alpha}(z).
\ee
The proof is complete.
\hfill$\Box$

\

{\it Proof of $(\mathrm{a}) \Rightarrow (\mathrm{d})$ for $p >1.$}\; Set $E = L^q(\mathbb{B}_n, \chi_{D(0, \gamma)} d \tau; \mathbb{C}^n).$ Consider the operator
\be
T (f)(z,w) =\big (\tilde{\nabla} f \big ) ( \varphi_z (w)), \quad f \in \mathcal{H} ( \mathbb{B}_n).
\ee
Note that $\varphi_z (D(0, \gamma)) = D(z, \gamma)$ and the measure $d \tau$ is invariant under any automorphism of $\mathbb{B}_n$ (cf. \cite[Proposition 1.13]{Zhu2005}), we have
\be\begin{split}
\| T (f) (z) \|_E & = \left ( \int_{\mathbb{B}_n} \big | \big ( \tilde{\nabla} f \big ) ( \varphi_z (w)) \big |^q \chi_{D(0, \gamma)} (w) d \tau (w) \right )^{\frac{1}{q}}\\
& = \left ( \int_{\mathbb{B}_n} \big | \tilde{\nabla} f ( w ) \big |^q \chi_{D(z, \gamma)} (w) d \tau (w) \right )^{\frac{1}{q}} = A^{\gamma, q}_{\tilde{\nabla}} (f) (z).
\end{split}\ee
On the other hand,
\be\begin{split}
A^{\gamma, q}_{\tilde{\nabla}} (f) (z) \le \left [ C_{\gamma}(1-|z|^2)^{-n-1} v(D(z, \gamma)) \right ]^{\frac{1}{2}} \|f\|_{\mathcal{B}} \le C \|f\|_{\mathcal{B}}.
\end{split}
\ee
This follows that $T$ is bounded from $\mathcal{B}$ into $L^{\8}_{\alpha} (\mathbb{B}_n, E).$ Notice that we have proved that $T$ is bounded from $\A^1_{\alpha}$ to $L^1_{\alpha} (\mathbb{B}_n, E).$ Thus, by Lemma \ref{le:InterpBloch} and the well known fact that
\be
L_p (\mathbb{B}_n, E) = (L^1_{\alpha} (\mathbb{B}_n, E), L^{\8}_{\alpha} (\mathbb{B}_n, E))_{\theta}\quad \text{with}\quad \theta = 1 - \frac{1}{p},
\ee
we conclude that $T$ is bounded from $\mathcal{A}^p_{\alpha}$ into $L^p_{\alpha} (\mathbb{B}_n, E)$ for any $1<p< \8,$ i.e.,
\be
\| A^{\gamma, q}_{\tilde{\nabla}} (f) \|_{p, \alpha} \le C \| f \|_{p,\alpha},\quad \forall f \in \mathcal{A}^p_{\alpha},
\ee
where $C$ depends only on $q, \g, n, p,$ and $\alpha.$ The proof is complete.
\hfill $\Box$

\begin{remark}\label{rk:p=8case}\rm
From the proofs of that $(b) \Longrightarrow (a)$ and that $(\mathrm{a}) \Rightarrow (\mathrm{d})$ for $p >1$ we find that Theorem \ref{th:AreaCharat} still holds true for the Bloch space. That is, for any $f \in \mathcal{H} (\mathbb{B}_n),$ $f \in \mathcal{B}$ if and only if one (or equivalently, all) of $A^{\gamma, q}_{\mathcal{R}} (f), A^{\gamma, q}_{\nabla}( f),$ and $A^{\gamma, q}_{\tilde{\nabla}} (f)$ is (or, are) in $L^{\8} (\mathbb{B}_n).$ Moreover,
\beq\label{eq:p=8AreaCharact}
\| f  \|_{\mathcal{B}} \approx \| A^{\gamma, q}_{\mathcal{R}} (f) \|_{L^{\8} (\mathbb{B}_n)} \approx \| A^{\gamma, q}_{\nabla}( f) \|_{L^{\8} (\mathbb{B}_n)} \approx \| A^{\gamma, q}_{\tilde{\nabla}} (f) \|_{L^{\8} (\mathbb{B}_n)},
\eeq
where ``$\approx$" depends only on $q, \gamma,$ and $n.$
\end{remark}

\section{Atomic decomposition for Bergman spaces}\label{Atomdecomp}

We let
$$
d(z,w)=\left | 1-\langle z,w\rangle \right |^{\frac{1}{2}},\; z , w \in \overline{\mathbb{B}}_n.
$$
It is known that $d$ satisfies the triangle inequality and the restriction of $d$ to ${\mathbb{S}_n}$ is a metric. As usual, $d$ is called the nonisotropic metric.

For any $\zeta\in{\mathbb{S}_n}$ and $r>0,$ the set
$$Q_r(\zeta)=\{z\in \mathbb{B}_n: d(z,\zeta)<r\}$$
is called a Carleson tube with respect to the nonisotropic metric $d.$ We usually write $Q = Q_r(\zeta)$ in short.

As usual, we define the atoms with respect to the Carleson tube as follows: for $1 <q<\infty,$ $a \in L^q (\mathbb{B}_n, d v_{\alpha})$ is said to be a $(1, q)_{\alpha}$-atom if there is a Carleson tube $Q$ such that
\begin{enumerate}[{\rm (1)}]

\item $a$ is supported in $Q;$

\item $\| a \|_{L^q (\mathbb{B}_n, d v_{\alpha})} \le v_{\alpha} (Q)^{\frac{1}{q}-1};$

\item $\int_{\mathbb{B}_n}a(z)\,dv_{\alpha}(z) = 0.$

\end{enumerate}
The constant function $1$ is also considered to be a $(1, q)_{\alpha}$-atom. Note that for any $(1, q)_{\alpha}$-atom $a,$
\be
\| a \|_{1, \alpha} = \int_Q | a | d v_{\alpha} \le v_{\alpha} (Q)^{1 - 1/q} \| a \|_{q, \alpha} \le 1.
\ee

Recall that $P_{\alpha}$ is the orthogonal projection from $L^2(\mathbb{B}_n,d v_{\alpha})$ onto ${\mathcal A}^2_{\alpha},$ which can be expressed as
\be
P_{\alpha}f(z)=\int_{\mathbb{B}_n}K^{\alpha}(z,w)f(w)dv_{\alpha}(w), \quad \forall f\in{L^1(\mathbb{B}_n,dv_{\alpha})},\alpha>-1,
\ee
where
\be
K^{\alpha}(z,w)=\frac{1}{\big(1-\langle z,w\rangle\big)^{n+1+\alpha}},\quad z,w\in\mathbb{B}_n.
\ee
$P_{\alpha}$ extends to a bounded projection from $L^p(\mathbb{B}_n,dv_{\alpha})$ onto $\mathcal A^p_{\alpha}$ ($1<p<\infty$).

We have the following useful estimates.

\begin{lemma}\label{le:AtEst}
For $\alpha> - 1$ and $1 < q < \8$ there exists a constant $C_{q, \alpha, n} >0$ such that
\be
\|P_{\alpha}(a)\|_{1,\,\alpha}\le C_{q, \alpha, n}
\ee
for any $(1, q)_{\alpha}$-atom $a.$
\end{lemma}

To prove Lemma \ref{le:AtEst}, we need first to show an inequality for reproducing kernel $K^{\alpha}$ associated with $d,$ which is essentially borrowed from \cite[Proposition 2.13]{Tch2008}.

\begin{lemma}\label{le:KerEst}
For $\alpha>-1$ there exists a constant $\delta >0$ such that for all $z,w\in{\mathbb{B}_n},\zeta\in{\mathbb{S}_n}$ satisfying $d(z,\zeta)> \delta d(w,\zeta),$ we have
\be
\left | K^{\alpha}(z,w)-K^{\alpha}(z,\zeta) \right |\le C_{\alpha, n} \frac{d(w,\zeta)}{d(z,\zeta)^{2(n+1+\alpha)+1}}.
\ee
\end{lemma}

\begin{proof} Note that
\be
K^{\alpha}(z,w)-K^{\alpha}(z,\zeta)=\int^1_0 \frac{d}{dt} \left ( \frac{1}{(1-\langle z,{\zeta}\rangle-t\langle z,{w-\zeta}\rangle)^{n+1+\alpha}} \right ) d t.
\ee
We have
\be
\left | K^{\alpha}(z,w)-K^{\alpha}(z,\zeta) \right | \le \int^1_0\frac{(n+1+\alpha)|\langle z,{w-\zeta}\rangle|}{|1-\langle z,{\zeta}\rangle-t\langle z,{w-\zeta}\rangle|^{n+2+\alpha}}dt.
\ee

Write $z=z_1+z_2$ and $w=w_1+w_2$, where $z_1$ and $w_1$ are parallel to $\zeta$, while $z_2$ and $w_2$ are perpendicular to $\zeta$. Then
$$\langle z,w\rangle-\langle z,{\zeta}\rangle=\langle z_2,w_2\rangle-\langle z_1,{w_1-\zeta}\rangle$$
and so
$$|\langle z,w\rangle-\langle z,{\zeta}\rangle|\le|z_2||w_2|+|w_1-\zeta|.$$
Since $|w_1-\zeta|=|1-\langle w,\zeta\rangle|,$
\be \begin{split}
|z_2|^2& =|z|^2-|z_1|^2<1-|z_1|^2<(1+|z_1|)(1+|z_1|) \\
& \le |1-\langle {z_1},{\zeta} \rangle | = 2|1-\langle {z},{\zeta}\rangle|,
\end{split} \ee
and similarly
\be
|w_2|^2\le2|1-\langle {w},{\zeta}\rangle|,
\ee
we have
\be\begin{split}
|\langle z,w\rangle-\langle z, \zeta\rangle| & \le 2|1-\langle z,\zeta\rangle|^{1/2}|1-\langle w,\zeta\rangle|^{1/2}+|1-\langle w,\zeta\rangle|\\
&= 2d(w,\zeta)[d(z,\zeta)+d(w,\zeta)] \le 2 \Big ( 1+\frac{1}{\delta} \Big ) \frac{1}{\delta}d^2(z,\zeta).
\end{split}\ee
This concludes that there is $\delta >1$ such that
\be
|\langle z,w-\zeta\rangle|<\frac{1}{2}|1-\langle z,\zeta\rangle|, \; \forall z, w \in \mathbb{B}_n,\; \zeta \in \mathbb{S}_n,
\ee
whenever $d(z,\zeta)> \delta d(w,\zeta).$
Then, we have
\be
|1-\langle z,\zeta\rangle-t\langle z, w-\zeta\rangle| > |1-\langle z,\zeta\rangle|-t|\langle z,\zeta-w\rangle| >\frac{1}{2}|1-\langle z,\zeta\rangle|.
\ee
Therefore,
\be \begin{split}
\left | K^{\alpha}(z,w)-K^{\alpha}(z,\zeta) \right | & \le \frac{2^{n+3+\alpha}(n+1+\alpha)(1+1/\delta) d(w,\zeta)d(z,\zeta)}{|1-\langle z,\zeta\rangle|^{n+2+\alpha}}\\
& \le C_{\alpha, n} \frac{d(w,\zeta)}{d(z,\zeta)^{2(n+1+\alpha)+1}}
\end{split}\ee
and the lemma is proved.
\end{proof}


{\it Proof of Lemma \ref{le:AtEst}.}\;
When $a$ is the constant function $1,$ the result is clear. Thus we may suppose $a$ is a $(1, q)_{\alpha}$-atom. Let $a$ be supported in a Carleson tuber $Q_r(\zeta)$ and $\delta r \le{\sqrt{2}},$ where $\delta$ is the constant in Lemma \ref{le:KerEst}. Since $P_{\alpha}$ is a bounded operator on $L^q(\mathbb{B}_n,dv_{\alpha}),$ we have
\be\begin{split}
\int_{Q_{\delta r}} | P_{\alpha}( a )| d v_{\alpha}(z) & \le v_{\alpha}(Q_{\delta r})^{1 - \frac{1}{q}} \|P_{\alpha}(a) \|_{q,\,\alpha} \le \| P_{\alpha} \|_{L^q} v_{\alpha}(Q_{\delta r})^{1- \frac{1}{q}} \|a \|_{q,\alpha} \le \| P_{\alpha} \|_{L^q}.
\end{split}\ee
Next, if $d(z,\zeta)> \delta r$ then
\be\begin{split}
\Big | \int_{\mathbb{B}_n} & \frac{a (w)}{(1-\langle z,w\rangle)^{n+1+\alpha}} d v_{\alpha}(w) \Big | \\
& = \left |\int_{Q_r (\zeta)} a (w) \left [ \frac{1}{(1-\langle z,w\rangle)^{n+1+\alpha}} - \frac{1}{(1-\langle z,\zeta\rangle)^{n+1+\alpha}} \right ] dv_{\alpha}(w) \right | \\
& \le C \int_{Q_{r}(\zeta)}|a (w)|\frac{d(w,\zeta)}{d(z,\zeta)^{2(n+1+\alpha)+1}}dv_{\alpha}(w)\\
& \le C r \int_{Q_{r}(\zeta)}|a (w)|dv_{\alpha}(w)\frac{1}{d(z,\zeta)^{2(n+1+\alpha)+1}} \le \frac{C r}{d(z,\zeta)^{2(n+1+\alpha)+1}}.
\end{split}
\ee
Then
\be\begin{split}
\int_{d(z,\,\zeta) > \delta r} &| P_{\alpha}(a) | d v_{\alpha}(z) \le C r \int_{d(z,\,\zeta)> \delta r} \frac{1}{d(z,\,\zeta)^{2(n+1+\alpha)+1}} d v_{\alpha}(z)\\
& = C r \sum_{k \ge 0} \int_{2^k \delta r < d(z,\,\zeta) \le 2^{k+1} \delta r } \frac{1}{d(z,\,\zeta)^{2(n+1+\alpha)+1}} d v_{\alpha}(z)\\
& \le C r \sum_{k \ge 0} \frac{v_{\alpha} (Q_{2^{k+1} \delta r})}{(2^k \delta r)^{2(n+1+\alpha)+1}} \le C r\sum_{k=0}^\infty\frac{(2^{k+1} \delta r)^{2(n+1+\alpha)}}{(2^k \delta r)^{2(n+1+\alpha)+1}}\le C,
\end{split}\ee
where we have used the fact that $v_{\alpha} (Q_r) \approx r^{2(n+1+\alpha)}$ in the third inequality (cf. \cite[Corollary 5.24]{Zhu2005}).
Thus, we get
\be
\int_{\mathbb{B}_n}|P_{\alpha}(a)|d v_{\alpha}(z)=\int_{Q_{\delta r}} |P_{\alpha}(a)| d v_{\alpha}(z) + \int_{d(z,\,\zeta)> \delta r}|P_{\alpha}(h)|d v_{\alpha}(z) \le C,
\ee
where $C$ depends only on $q, n$ and $\alpha.$
\hfill$\Box$

\

Now we turn to the atomic decomposition of $\mathcal{A}^1_{\alpha}$ ($\alpha > -1$) with respect to the Carleson tubes. Recall that $\| a \|_{1, \alpha} \le 1$ for any $(1, q)_{\alpha}$-atom $a.$ Then, we define $\mathcal{A}^{1,q}_{\alpha}$ as the space of all $f \in \mathcal{A}^1_{\alpha}$ which admits a decomposition
\be
f=\sum_i \lambda_i P_{\alpha} a_i\quad \text{and}\quad \sum_i |\lambda_i | \le C_q \|f\|_{1,\,\alpha},
\ee
where for each $i,$ $a_i$ is an $(1,q)_{\alpha}$-atom and $\lambda_i \in \mathbb{C}$ so that $\sum_i | \lambda_i | < \8.$ We equip this space with the norm
\be
\|f\|_{\mathcal{A}^{1,q}_{\alpha}} = \inf \Big \{ \sum_i |\lambda_i |:\; f=\sum_i \lambda_i P_{\alpha} a_i \Big \}
\ee
where the infimum is taken over all decompositions of $f$ described above.

It is easy to see that  $\mathcal{A}^{1,q}_{\alpha}$ is a Banach space.
By Lemma \ref{le:AtEst} we have the contractive inclusion $\mathcal{A}^{1,q}_{\alpha} \subset \mathcal{A}^1_{\alpha}.$
We will prove in what follows that these two spaces coincide. That establishes the ``real-variable" atomic decomposition of the Bergman space $\mathcal{ A}^1_{\alpha}.$ In fact, we will show the remaining inclusion $\mathcal{A}^1_{\alpha} \subset \mathcal{A}^{1,q}_{\alpha}$ by duality.

\begin{theorem}\label{th:AtDecomp}
Let $1< q<\infty$ and $\alpha>-1.$ For every $f \in{\mathcal A^1_{\alpha}}$ there exist a sequence $\{a_i\}$ of $(1, q)_{\alpha}$-atoms and a sequence $\{\lambda_i\}$ of complex numbers such that
\beq \label{eq:AtDecomp}
f=\sum_i \lambda_i P_{\alpha} a_i\quad \text{and}\quad \sum_i |\lambda_i | \le C_q \|f\|_{1,\,\alpha}.
\eeq
Moreover,
\be
\|f\|_{1,\,\alpha} \approx \inf \sum_i |\lambda_i |
\ee
where the infimum is taken over all decompositions of $f$ described above and $`` \approx "$ depends only on $\alpha, n,$ and $q.$
\end{theorem}

Recall that the dual space of $\mathcal{A}^1_{\alpha}$ is the Bloch space $\mathcal{B}$ (we refer to \cite {Zhu2005} for details). The Banach dual of $\mathcal A^1_{\alpha}$ can be identified with $\mathcal B$ (with equivalent norms) under the integral pairing
\be
\langle f,g \rangle_{\alpha}= \lim_{r \to 1^-}\int_{\mathbb{B}_n} f(r z) \overline{g(z)} d v_{\alpha}(z),\quad f \in{\mathcal A^1_{\alpha}},\; g \in \mathcal{B}.
\ee
(cf. \cite[Theorem 3.17]{Zhu2005}.)

In order to prove Theorem \ref{th:AtDecomp}, we need the following result, which can be found in \cite{BBCZ1990} (see also \cite[Theorem 5.25]{Zhu2005}).

\begin{lemma}\label{le:Bloch-BMO}
Suppose $\alpha>-1$ and $1 \le p < \8.$ Then, for any $f \in \mathcal{H}(\mathbb{B}_n),$ $f$ is in $\mathcal B$ if and only if there exists a constant $C>0$ depending only on $\alpha$ and $p$ such that
\be
\frac{1}{v_{\alpha}(Q_{r}(\zeta))}\int_{Q_{r}(\zeta)}|f-f_{\alpha,\,Q_{r}(\zeta)}|^p d v_{\alpha}\le C
\ee
for all $r>0$ and all $\zeta\in\mathbb{S}_n,$ where
$$f_{\alpha,\,Q_{r}(\zeta)}=\frac{1}{Q_{r}(\zeta)}\int_{Q_{r}(\zeta)}f(z)dv_{\alpha}(z).$$
Moreover,
\be
\| f \|_{\mathcal{B}} \approx \sup_{r > 0, \zeta \in \mathbb{S}} \Big ( \frac{1}{v_{\alpha}(Q_{r}(\zeta))}\int_{Q_{r}(\zeta)}|f-f_{\alpha,\,Q_{r}(\zeta)}|^p d v_{\alpha} \Big )^{\frac{1}{p}},
\ee
where ``$\approx$" depends only on $\alpha, p,$ and $n.$
\end{lemma}

As noted above, we will prove Theorem \ref{th:AtDecomp} via duality. To this end, we first prove the following duality theorem.

\begin{proposition}\label{prop:AtDuality}
For any $1 < q < \8$ and $\alpha > -1,$ we have $(\mathcal{A}^{1,q}_{\alpha})^* = \mathcal{B}$ isometrically. More precisely,
\begin{enumerate}[{\rm (i)}]

\item Every $g \in \mathcal{B}$ defines a continuous linear functional $\varphi_g$ on $\mathcal{A}^{1,q}_{\alpha}$ by
\beq\label{eq:duality_bracket_at}
\varphi_g (f) = \lim_{r \to 1^-}\int_{\mathbb{B}_n} f (r z) \overline{g(z)} d v_{\alpha} (z),\quad \forall f \in \mathcal{A}^{1,q}_{\alpha}.
\eeq

\item Conversely, each $\varphi \in (\mathcal{A}^{1,q}_{\alpha})^*$ is given as \eqref{eq:duality_bracket_at}
by some $g \in \mathcal{B}.$

\end{enumerate}
Moreover, we have
\beq\label{eq:BlochDualityNorm}
\| \varphi_g \| \approx |g(0)| + \| g \|_{\mathcal{B}},\quad \forall g \in \mathcal{B}.
\eeq
\end{proposition}

\begin{proof}
Let $p$ be the conjugate index of $q,$ i.e., $1/p+1/q=1.$ We first show $\mathcal{B} \subset (\mathcal{A}^{1,q}_{\alpha})^*.$ Let $g \in \mathcal{B}.$ For any $(1,q)_{\alpha}$-atom $a,$ by Lemma \ref{le:Bloch-BMO} we have
\be\begin{split}
\left |\int_{\mathbb{B}_n} P_{\alpha} a (z) \overline{g(z)} dv_{\alpha} (z) \right |& = | \langle P_{\alpha}(a_j),g \rangle_{\alpha} | = \left | \int_{\mathbb{B}_n} a \overline g dv_{\alpha} \right| = \left|\int_{\mathbb{B}_n} a \overline{\left( g-g_Q \right)}dv_{\alpha}\right|\\
& \le \left(\int_Q |a|^q dv_{\alpha}\right)^{1/q} \left ( \int_Q |g-g_Q|^p d v_{\alpha} \right)^{1/p}\\
& \le \left( \frac{1}{v_{\alpha}(Q)} \int_Q |g-g_Q|^p d v_{\alpha} \right)^{1/p} \le C \| g \|_{\mathcal B}.
\end{split}\ee
On the other hand, for the constant function $1$ we have $P_{\alpha} 1 = 1$ and so
\be
\left |\int_{\mathbb{B}_n} P_{\alpha} 1 (z) \overline{g(z)} dv_{\alpha} (z) \right | = \left |\int_{\mathbb{B}_n} g(z) dv_{\alpha} (z) \right | = | g(0)|.
\ee
Thus, we deduce that
\be
\left|\int_{\mathbb{B}_n} f \bar{g} d v_{\alpha}\right|\le C \|f\|_{\mathcal{A}_{\alpha}^{1,q}} (|g(0)| +  \|g\|_{\mathcal B})
\ee
for any finite linear combination $f$ of $(1,q)_{\alpha}$-atoms. Hence, $g$ defines a continuous linear functional $\varphi_g$ on a dense
subspace of ${\mathcal{A}_{\alpha}^{1,\,q}}$ and $\varphi_g$ extends to a continuous linear functional on ${\mathcal{A}_{\alpha}^{1,\,q}}$ such that
\be
\left|\varphi_g (f) \right | \le C (|g(0)|+ \|g\|_{\mathcal B}) \| f \|_{{\mathcal{A}_{\alpha}^{1,\,q}}}
\ee
for all $f \in \mathcal{A}_{\alpha}^{1,\,q}.$

Next let $\varphi$ be a bounded linear functional on $\mathcal{A}_{\alpha}^{1,\,q}.$ Note that
\be
\mathcal{H}^q (\mathbb{B}_n, d v_{\alpha}) = \mathcal{H} (\mathbb{B}_n) \cap L^q (\mathbb{B}_n, d v_{\alpha}) \subset \mathcal{A}_{\alpha}^{1,\,q}.
\ee
Then, $\varphi$ is a bounded linear functional on $\mathcal{H}^q (\mathbb{B}_n, d v_{\alpha}).$ By duality there exists $g \in \mathcal{H}^p (\mathbb{B}_n, d v_{\alpha})$ such that
\be
\varphi (f) = \int_{\mathbb{B}_n} f \bar{g} d v_{\alpha}, \quad \forall f \in \mathcal{H}^q (\mathbb{B}_n, d v_{\alpha}).
\ee
Let $Q=Q_{r}(\zeta)$ be a Carleson tube. For any $f \in L^q (\mathbb{B}_n, d v_{\alpha})$ supported in $Q,$ it is easy to check that
\be
a_f = (f-f_Q) \chi_Q / [\| f \|_{L^q} v_{\alpha} (Q)^{1/p} ]
\ee
is a $(1,q)$-atom. Then, $| \varphi (P_{\alpha} a_f) | \le \| \varphi \|$ and so
\be
\big | \varphi (P_{\alpha} [(f-f_Q) \chi_Q]) \big | \le \| \varphi \| \| f \|_{L^q} v_{\alpha} (Q)^{1/p}.
\ee
Hence, for any $f \in L^q (\mathbb{B}_n, d v_{\alpha})$ we have
\be\begin{split}
\left | \int_Q f \overline{(g-g_Q)} d v_{\alpha} \right | & = \left | \int_Q (f-f_Q) \bar{g} d v_{\alpha} \right |  = \left | \int_{\mathbb{B}_n} (f-f_Q) \chi_Q \bar{g} dv_{\alpha} \right |\\
&  = \left | \int_{\mathbb{B}_n} P_{\alpha}[ (f-f_{Q})\chi_Q] \bar{g} d v_{\alpha} \right | = \left |\varphi ( P_{\alpha}[ (f-f_{Q}) \chi_{Q}]) \right |\\
& \le \|\varphi \| \left \|(f-f_{Q}) \chi_{Q} \right \|_{L^q({\mathbb{B}_n},\,dv_{\alpha})}{v_{\alpha}(Q)}^{1/p} \le 2 \| \varphi\| \left \|f \right \|_{L^q(Q,\,dv_{\alpha})}{v_{\alpha}(Q)}^{1/p}.
\end{split}\ee
This concludes that
\be
\left(\frac{1}{v_{\alpha}(Q)} \int_Q \left |g - g_Q \right |^p d v_{\alpha} \right)^{1/p} \le 2 \|\varphi \|.
\ee
By Lemma \ref{le:Bloch-BMO} we have that $g \in \mathcal B$ and $\|g\|_{\mathcal B}\le C\|\varphi \|.$ Therefore, $\varphi$ is given as \eqref{eq:duality_bracket_at} by $g$ with $|g(0)| + \| g \|_{\mathcal{B}} \le C \| \varphi \|.$
\end{proof}

Now we are ready to prove Theorem \ref{th:AtDecomp}.

\

{\it Proof of Theorem \ref{th:AtDecomp}.}\; By Lemma \ref{le:AtEst} we know that $\mathcal{A}^{1,q}_{\alpha} \subset \mathcal{A}^1_{\alpha}.$ On the other hand, by Proposition \ref{prop:AtDuality} we have $(\mathcal{A}^1_{\alpha})^* = (\mathcal{A}^{1,q}_{\alpha})^*.$ Hence, by duality we have $\| f \|_{1,q} \approx \| f \|_{\mathcal{A}^{1,q}_{\alpha}}.$ \hfill$\Box$

\begin{remark}\label{rk:atomdecomp-p<1}\rm
\begin{enumerate}

\item One would like to expect that when $0 < p < 1,$ $\A^p_{\alpha}$ also admits an atomic decomposition in terms of atoms with respect to Carleson tubes. However, the proof of Theorem \ref{th:AtDecomp} via duality cannot be extended to the case $0< p<1.$ At the time of this writing, this problem is entirely open.

\item The real-variable atomic decomposition of Bergman spaces should be known to specialists in the case $p=1.$ Indeed, based on their theory of harmonic analysis on homogeneous spaces, Coifman and Weiss \cite{CW1977} claimed that the Bergman space $\A^1$ admits an atomic decomposition in terms of atoms with respect to $(\mathbb{B}_n, \varrho, d v),$ where
\be
\varrho (z, w) = \left \{\begin{split}
& \big | |z| - |w| \big | + \Big | 1 - \frac{1}{|z| |w|} \langle z, w \rangle \Big |,\quad \text{if}\; z, w \in \mathbb{B}_n \backslash \{0\},\\
& |z| + |w|,\quad \text{otherwise}.
\end{split} \right.
\ee
This also applies to $\A^1_{\alpha}$ because $(\mathbb{B}_n, \varrho, d v_{\alpha})$ is a homogeneous space for $\alpha > -1$ (see e.g. \cite{Tch2008}). However, the approach of Coifman and Weiss is again based on duality and therefore not constructive and cannot be applied to the case $0<p<1.$ Recently, the present authors \cite{CO2013} extend this result to the case $0<p <1$ through using a constructive method.

\end{enumerate}
\end{remark}

\section{Area integral inequalities: Real-variable methods}\label{AreaCharactp=1}

In this section, we will prove the area integral inequality for the Bergman space $\A^1_{\alpha}$ via atomic decomposition established in Section \ref{Atomdecomp}.

\begin{theorem}\label{th:AreaCharactp=1}
Suppose $1< q < \8, \gamma > 0,$ and $\alpha > -1.$ Then,
\beq\label{eq:AreaCharactp=1}
\| A_{\tilde{\nabla}}^{\gamma, q} (f) \|_{1, \alpha} \lesssim \| f \|_{1, \alpha},\quad \forall f \in \mathcal{H} (\mathbb{B}_n).
\eeq
\end{theorem}

This is the assertion $(\mathrm{a}) \Rightarrow (\mathrm{d})$ of Theorem \ref{th:AreaCharat} in the case $p=1.$ The novelty of the proof here is to involve a real-variable method.

The following lemma is elementary.

\begin{lemma}\label{le:L2InGraEst}
Suppose $1 < q < \8, \gamma>0$ and $\alpha > -1.$ If $f \in \mathcal{A}^q_{\alpha},$ then
\be
\int_{\mathbb{B}_n} | A^{\gamma, q}_{\tilde{\nabla}} (f) (z) |^q d v_{\alpha} \approx \int_{\mathbb{B}_n} | f (z) - f(0) |^q d v_{\alpha},
\ee
where $``\approx "$ depends only on $q, \gamma, \alpha, $ and $n.$
\end{lemma}

\begin{proof}
Note that $v_{\alpha} (D(z,\gamma)) \approx (1-|z|^2)^{n+1+ \alpha}.$ Then
\be\begin{split}
\int_{\mathbb{B}_n} | A^{\gamma, q}_{\tilde{\nabla}} (f) (z) |^q d v_{\alpha} & = \int_{\mathbb{B}_n} \int_{D(z,\gamma)} (1-|w|^2)^{-1-n} | \tilde{\nabla} f (w) |^q d v (w) d v_{\alpha}(z)\\
& = \int_{\mathbb{B}_n} v_{\alpha} (D(w, \gamma)) (1-|w|^2)^{-1-n} | \tilde{\nabla} f (w) |^q d v (w)\\
& \approx \int_{\mathbb{B}_n} | \tilde{\nabla} f (w) |^q d v_{\alpha} (w) \approx \int_{\mathbb{B}_n} | f(w) - f(0)|^q d v_{\alpha} (w).
\end{split}\ee
In the last step we have used \cite[Theorem 2.16 (b)]{Zhu2005}.
\end{proof}

\

{\it Proof of \eqref{eq:AreaCharactp=1}}.\; By Theorem \ref{th:AtDecomp}, it suffices to show that for $1< q < \8, \gamma>0,$ and $ \alpha>-1$ there exists $C>0$ such that
\be
\| A^{\gamma, q}_{\tilde{\nabla}} (P_{\alpha} a) \|_{1, \alpha} \le C
\ee
for all $(1,q)_{\alpha}$-atoms $a.$ Given an $(1,q)_{\alpha}$-atom $a$ supported in $Q = Q_r (\zeta),$ by Lemma \ref{le:L2InGraEst} we have
\be\begin{split}
\int_{2 Q} A^{\gamma, q}_{\tilde{\nabla}} (P_{\alpha} a) d v_{\alpha} & \le v_{\alpha}( 2 Q)^{1-\frac{1}{q}} \Big ( \int_{2 Q} \big [ A^{\gamma, q}_{\tilde{\nabla}} (P_{\alpha} a) \big ]^q d v_{\alpha} \Big )^{\frac{1}{q}}\\
& \le C v_{\alpha}(Q)^{1-\frac{1}{q}} \Big ( \int_{\mathbb{B}_n} | P_{\alpha} a (z) - P_{\alpha} a (0) |^q d v_{\alpha} \Big )^{\frac{1}{q}}\\
& \le C v_{\alpha}(Q)^{1-\frac{1}{q}} \| a \|_{q, \alpha} \le C,
\end{split}\ee
where $2 Q = Q_{2 r} (\zeta).$ On the other hand,
\be\begin{split}
\int_{(2 Q)^c} & A^{\gamma, q}_{\tilde{\nabla}} (P_{\alpha} a) d v_{\alpha}
= \int_{(2 Q)^c} \Big ( \int_{D(z, \gamma)}  |\tilde{\nabla} P_{\alpha} a(w)|^q d \tau (w) \Big )^{\frac{1}{q}} d v_{\alpha} (z)\\
= & \int_{(2 Q)^c} \Big ( \int_{D(z, \gamma)} \Big | \int_Q \tilde{\nabla}_w [ K^{\alpha} (w,u) - K^{\alpha} (w, \zeta)] a(u) d v_{\alpha}(u) \Big |^q d \tau (w) \Big )^{\frac{1}{q}} d v_{\alpha} (z)\\
\le & \| a \|_{q, \alpha} \int_{(2 Q)^c} \Big ( \int_{D(z, \gamma)} \Big ( \int_Q |\tilde{\nabla}_w [ K^{\alpha} (w,u) - K^{\alpha} (w, \zeta)]|^{\frac{q}{q-1}} d v_{\alpha}(u) \Big )^{q-1} d \tau (w) \Big )^{\frac{1}{q}} d v_{\alpha} (z)\\
\le & \int_{(2 Q)^c} \Big ( \int_{D(z, \gamma)} \sup_{u \in Q} |\tilde{\nabla}_w [ K^{\alpha} (w,u) - K^{\alpha} (w, \zeta)]|^q d \tau (w) \Big )^{\frac{1}{q}} d v_{\alpha} (z),
\end{split}\ee
where $(2 Q)^c = \mathbb{B}_n \backslash 2 Q.$

An immediate computation yields that
\be\begin{split}
{\nabla}_{w}& [K^{\alpha}(w,u)-K^{\alpha}(w,\zeta)]\\
=& (n+1+\alpha)\Big [ \frac{\bar{u}}{(1-\langle w,u\rangle)^{n+2+\alpha}}-\frac{\bar{\zeta}}{(1-\langle w,\zeta \rangle)^{n+2+\alpha}} \Big ]\\
=& (n+1+\alpha)\frac{\bar{u}(1-\langle w,\zeta \rangle)^{n+2+\alpha}-\bar{\zeta}(1-\langle w,u\rangle)^{n+2+\alpha}}{(1-\langle
w,u\rangle)^{n+2+\alpha}(1-\langle w,\zeta \rangle)^{n+2+\alpha}}
\end{split}\ee
and
\be\begin{split}
\mathcal{R}_{w} & [K^{\alpha}(w,u)-K^{\alpha}(w,\zeta)]\\
=&(n+1+\alpha) \Big [ \frac{\langle w,u\rangle}{(1-\langle w,u\rangle)^{n+2+\alpha}}-\frac{\langle w,\zeta \rangle}{(1-\langle w,\zeta \rangle)^{n+2+\alpha}} \Big ]\\
=&(n+1+\alpha) \frac{\langle w,u\rangle(1-\langle w, \zeta \rangle)^{n+2+\alpha}-\langle w,\zeta \rangle (1-\langle w,u\rangle)^{n+2+\alpha}}{(1-\langle
w,u\rangle)^{n+2+\alpha}(1-\langle w,\zeta \rangle)^{n+2+\alpha}}.
\end{split}\ee
Moreover,
\be\begin{split}
\big | {\nabla}_{w} & [K^{\alpha}(w,u)-K^{\alpha}(w,\zeta)] \big |^{2}\\
= &(n+1+\alpha)^{2} \Big \{ \frac{|u|^2 |1-\langle w,\zeta \rangle|^{2(n+2+\alpha)} + |1-\langle w,u\rangle|^{2(n+2+\alpha)}}{|1-\langle
w,u\rangle|^{2(n+2+\alpha)}|1-\langle w,\zeta \rangle|^{2(n+2+\alpha)}}\\
&\quad - \frac{(1-\langle w,\zeta\rangle)^{n+2+\alpha} ( 1- \langle
u, w \rangle)^{n+2+\alpha} \langle \zeta,u\rangle}{|1-\langle w,u\rangle|^{2(n+2+\alpha)}|1-\langle w,\zeta
\rangle|^{2(n+2+\alpha)}}\\
&\quad -\frac{(1-\langle w,u\rangle)^{n+2+\alpha} (1-\langle \zeta, w \rangle)^{n+2+\alpha} \langle u,\zeta\rangle}{|1-\langle
w,u\rangle|^{2(n+2+\alpha)}|1-\langle w,\zeta \rangle|^{2(n+2+\alpha)}} \Big \},
\end{split}\ee
and
\be\begin{split}
& \big | \mathcal{R}_{w}[K^{\alpha}(w,u)-K^{\alpha}(w,\zeta)] \big |^2\\
& = (n+1+\alpha)^{2} \Big \{ \frac{|\langle w,u\rangle|^{2}|1-\langle w,\zeta \rangle|^{2(n+2+\alpha)}+|\langle w,\zeta \rangle|^{2}|1-\langle
w,u\rangle|^{2(n+2+\alpha)}}{|1-\langle w,u\rangle|^{2(n+2+\alpha)}|1-\langle w,\zeta
\rangle|^{2(n+2+\alpha)}}\\
& \quad - \frac{\langle w,u\rangle \langle \zeta, w \rangle (1-\langle w,\zeta\rangle )^{n+2+\alpha} (1-\langle
u, w \rangle )^{n+2+\alpha}}{|1-\langle w,u\rangle|^{2(n+2+\alpha)}|1-\langle w,\zeta \rangle|^{2(n+2+\alpha)}}\\
& \quad - \frac{\langle w,\zeta\rangle \langle u, w \rangle (1-\langle w,u \rangle )^{n+2+\alpha} (1-\langle
\zeta, w \rangle)^{n+2+\alpha}}{|1-\langle w,u\rangle|^{2(n+2+\alpha)}|1-\langle w,\zeta \rangle|^{2(n+2+\alpha)}} \Big \}.
\end{split}\ee
Then, we have
\be\begin{split}
\big | {\nabla}_{w} & [K^{\alpha}(w,u)-K^{\alpha}(w,\zeta)] \big |^2 - \big | \mathcal{R}_{w}[K^{\alpha}(w,u)-K^{\alpha}(w,\zeta)] \big |^2\\
& = \frac{(n+1+\alpha)^2}{|1-\langle w,u\rangle|^{2(n+2+\alpha)} |1-\langle w,\zeta \rangle|^{2(n+2+\alpha)}}\\
& \quad \times \Big \{ (|u|^{2}-|\langle w,u\rangle|^2) |1 - \langle w,\zeta\rangle|^{2(n+2+\alpha)}\\
& \quad + (1-|\langle w,\zeta\rangle|^{2})|1-\langle w,u\rangle|^{2(n+2+\alpha)} \\
& \quad + (\langle w,u \rangle \langle \zeta, w\rangle - \langle \zeta,u\rangle)(1-\langle w,\zeta\rangle)^{n+2+\alpha} (1-\langle
u, w \rangle)^{n+2+\alpha} \\
& \quad + (\langle w,\zeta\rangle \langle u, w \rangle - \langle u,\zeta \rangle)(1-\langle
w,u \rangle)^{n+2+\alpha} (1-\langle \zeta, w \rangle)^{n+2+\alpha} \Big \}.
\end{split}\ee
Note that for any $f \in \mathcal{H} (\mathbb{B}_n),$
\be
| \tilde{\nabla} f (z) |^2 = (1-|z|^2)(| \nabla f (z) |^2 - | \mathcal{R} f (z)|^2 ), \quad z \in \mathbb{B}_n
\ee
(cf. \cite[Lemma 2.13]{Zhu2005}). It is concluded that
\be\begin{split}
\big | \tilde{\nabla}_w & [K^{\alpha}(w,u)-K^{\alpha}(w,\zeta)] \big |\\
& \le \frac{(n+1+\alpha) (1- |w|^2)^{\frac{1}{2}}}{|1-\langle w,u\rangle|^{n+2+\alpha} |1-\langle w,\zeta \rangle|^{n+2+\alpha}}\\
& \quad \times \Big \{ (1 -|\langle w,u\rangle|^2) |1 - \langle w,\zeta\rangle|^{2(n+2+\alpha)} + (1-|\langle w,\zeta\rangle|^{2})|1-\langle w,u\rangle|^{2(n+2+\alpha)} \\
& \quad + \big [ \langle w,u-\zeta \rangle \langle \zeta, w\rangle + (|\langle w, \zeta \rangle|^2- 1) + (1- \langle \zeta,u\rangle) \big ]\\
& \quad \;\; \times (1-\langle w,\zeta\rangle)^{n+2+\alpha} (1-\langle u, w \rangle)^{n+2+\alpha} \\
& \quad + \big [ \langle w,\zeta -u \rangle \langle u, w \rangle + ( | \langle w, u \rangle |^2 - 1) + (1 - \langle u,\zeta \rangle) \big ]\\
& \quad \; \; \times (1-\langle w,u \rangle)^{n+2+\alpha} (1-\langle \zeta, w \rangle)^{n+2+\alpha} \Big \}^{\frac{1}{2}}\\
& \le \frac{(n+1+\alpha) (1-|w|^2)^{\frac{1}{2}}(M_1 + M_2 + M_3 + M_4)^{\frac{1}{2}}}{|1-\langle w,u\rangle|^{n+2+\alpha} |1-\langle w,\zeta \rangle|^{ n+2+\alpha}},
\end{split}\ee
where
\be\begin{split}
M_1 =&  |1-\langle w,\zeta\rangle |^{n+2+\alpha} | 1-\langle u, w \rangle |^{n+2+\alpha} \left | \langle w,u-\zeta \rangle \langle \zeta, w \rangle +(1-\langle \zeta,u\rangle) \right |,\\
M_2 =& | 1-\langle w,u \rangle |^{n+2 + \alpha} | 1-\langle \zeta, w \rangle |^{n+2+\alpha} \left | \langle w,\zeta - u \rangle \langle u, w \rangle + (1-\langle u,\zeta \rangle) \right |,\\
M_3 =&  (1- |\langle w,u\rangle|^2) |1-\langle \zeta, w \rangle |^{n+2+\alpha} \left | (1-\langle w, \zeta \rangle)^{n+2+\alpha}- (1-\langle w,u \rangle )^{n+2+\alpha} \right |,\\
M_4 =&  (1-|\langle w,\zeta\rangle|^2) |1-\langle u, w \rangle |^{n+2+\alpha} \left | (1-\langle w,u \rangle)^{n+2+\alpha} - (1- \langle w,\zeta\rangle)^{n+2+\alpha} \right |,
\end{split}\ee
for $w \in D(z, \g),\;u \in Q_r (\zeta)$ and $z \in \mathbb{B}_n, \zeta \in \mathbb{S}_n.$

Hence,
\be\begin{split}
\int_{(2 Q)^{c}} A^{\gamma, q}_{\tilde{\nabla} } (P_{\alpha} a) d v_{\alpha}
\leq & \int_{(2 Q)^{c}} \Big ( \int_{D(z,\gamma)} \sup_{u \in Q} \big |\tilde{\nabla}_{w} [K^{\alpha}(w,u)-K^{\alpha}(w,\zeta) ]\big |^q d \tau (w) \Big )^{\frac{1}{q}} dv_{\alpha}(z)\\
\leq & (n+1+\alpha) \int_{(2 Q)^{c}} (I_1+I_2+I_3+I_4) dv_{\alpha} (z),
\end{split}\ee
where
\be\begin{split}
I_1 & = \Big ( \int_{D(z,\gamma)}\sup_{u \in Q}\frac{(1-|w|^{2})^{\frac{q}{2}} M_1^{\frac{q}{2}}}{|1-\langle
w,u\rangle|^{q(n+2+\alpha)}|1-\langle w,\zeta \rangle|^{q(n+2+\alpha)}} d \tau (w) \Big )^{\frac{1}{q}},\\
I_2 & = \Big ( \int_{D(z,\gamma)}\sup_{u \in Q} \frac{(1-|w|^{2})^{\frac{q}{2}} M_2^{\frac{q}{2}}}{|1-\langle
w,u\rangle|^{q(n+2+\alpha)}|1-\langle w,\zeta \rangle|^{q(n+2+\alpha)}}d \tau (w) \Big )^{\frac{1}{q}},\\
I_3 & = \Big ( \int_{D(z,\gamma)}\sup_{u \in Q} \frac{(1-|w|^2)^{\frac{q}{2}} M_3^{\frac{q}{2}}}{|1-\langle
w,u\rangle|^{q(n+2+\alpha)}|1-\langle w,\zeta \rangle|^{q(n+2+\alpha)}} d \tau (w) \Big )^{\frac{1}{q}},\\
I_4 & = \Big ( \int_{D(z,\gamma)}\sup_{u \in Q} \frac{(1-|w|^2)^{\frac{q}{2}} M_4^{\frac{q}{2}}}{|1-\langle w,u\rangle|^{q(n+2+\alpha)}|1-\langle w,\zeta
\rangle|^{q(n+2+\alpha)}} d \tau (w) \Big )^{\frac{1}{q}}.
\end{split}\ee

We first estimate $I_1.$ Note that
\be\begin{split}
M_1 &\leq (|\langle w,u-\zeta\rangle| + |1-\langle \zeta,u\rangle|) |1-\langle w,\zeta\rangle|^{n+2+\alpha}|1-\langle w,u\rangle|^{n+2+\alpha}\\
& \leq \left ( 2|1-\langle u,\zeta\rangle|^{\frac{1}{2}} \big ( |1-\langle w,\zeta\rangle|^{\frac{1}{2}}+|1-\langle
u,\zeta\rangle|^{\frac{1}{2}} \big )+|1-\langle \zeta,u\rangle| \right )\\
& \quad \times |1-\langle w,\zeta\rangle|^{n+2+\alpha}|1-\langle w,u\rangle|^{n+2+\alpha}\\
& \leq \left ( 2|1-\langle u,\zeta\rangle|^{\frac{1}{2}} \big ( C_\gamma|1-\langle z,\zeta\rangle|^{\frac{1}{2}}+\frac{1}{2}|1-\langle
z,\zeta\rangle|^{\frac{1}{2}} \big ) + |1-\langle \zeta,u\rangle| \right )\\
& \quad \times |1-\langle w,\zeta\rangle|^{n+2+\alpha}|1-\langle w,u\rangle|^{n+2+\alpha}\\
& \leq \left ( C_{\g} r |1-\langle z,\zeta\rangle|^{\frac{1}{2}}+r^2 \right ) |1-\langle
w,\zeta\rangle|^{n+2+\alpha}|1-\langle w,u\rangle|^{n+2+\alpha},
\end{split}\ee
where the second inequality is the consequence of the following fact
which has appeared in the proof of Lemma \ref{le:KerEst}
\be
|\langle w,u-\zeta\rangle|\leq2|1-\langle u,\zeta\rangle|^{\frac{1}{2}}(|1-\langle
w,\zeta\rangle|^{\frac{1}{2}}+|1-\langle
u,\zeta\rangle|^{\frac{1}{2}});
\ee
the third inequality is obtained by Lemma \ref{le:KerEstimation2} and the fact
\be
|1-\langle u,\zeta\rangle|^{\frac{1}{2}}<r<\frac{1}{2}|1-\langle
z,\zeta\rangle|^{\frac{1}{2}}
\ee
for $u\in Q$ and $z\in(2 Q)^{c}$. Since
\be\begin{split}
|1-\langle z,u\rangle|^{\frac{1}{2}} & \ge|1-\langle z,\zeta\rangle|^{\frac{1}{2}}-|1-\langle u,\zeta\rangle|^{\frac{1}{2}}\\
& \ge |1-\langle z,\zeta\rangle|^{\frac{1}{2}}-\frac{1}{2}|1-\langle z,\zeta\rangle|^{\frac{1}{2}} \ge \frac{1}{2}|1-\langle z,\zeta\rangle|^{\frac{1}{2}},
\end{split}\ee
by Lemmas \ref{le:KerEstimation} and \ref{le:KerEstimation2} we have
\be\begin{split}
I_{1}\leq & \left ( \int_{D(z,\gamma)}\sup_{u \in Q} \frac{(1-|w|^{2})^{\frac{q}{2}} \big [ C r |1-\langle
z,\zeta\rangle|^{\frac{1}{2}}+ r^{2} \big ]^{\frac{q}{2}}}{|1-\langle w,u\rangle|^{\frac{q}{2}(n+2+\alpha)}|1-\langle w,\zeta
\rangle|^{\frac{q}{2} (n+2+\alpha)}}d \tau (w) \right)^{\frac{1}{q}}\\
\leq & C_{\g} \left ( \int_{D(z,\gamma)} \sup_{u \in Q} \frac{(1-|z|^{2})^{\frac{q}{2}} \big ( C r |1-\langle
z,\zeta\rangle|^{\frac{1}{2}}+r^{2} \big )^{\frac{q}{2}}}{|1-\langle z,u\rangle|^{\frac{q}{2}( n+2+\alpha)}|1-\langle z,\zeta
\rangle|^{\frac{q}{2} (n+2+\alpha )}}d \tau (w) \right)^{\frac{1}{q}}\\
\leq & C_{\g} \left ( \frac{(1-|z|^{2})^{\frac{q}{2}} \big ( r |1-\langle
z,\zeta\rangle|^{\frac{1}{2}}+r^{2} \big )^{\frac{q}{2}}}{|1-\langle z,\zeta \rangle|^{q(n+2+\alpha)}} \right)^{\frac{1}{q}}
\le C_{\g} \frac{ \big ( r|1-\langle z,\zeta\rangle|^{\frac{1}{2}}+r^{2} \big )^{\frac{1}{2}}}{|1-\langle z,\zeta \rangle|^{n+ \frac{3}{2}+\alpha}}\\
\leq & C_{\g} \Big ( \frac{r^{\frac{1}{2}}}{d (z, \zeta)^{2(n+ 1 +\alpha )+ \frac{1}{2}}} + \frac{r}{d (z, \zeta)^{2(n+ 1 +\alpha)+1}} \Big ).\\
\end{split}\ee
Hence,
\be\begin{split}
\int_{(2 Q)^{c}}I_1 dv_{\alpha}(z) \lesssim & \int_{(2 Q)^{c}}\frac{r^{\frac{1}{2}}}{d(z,\zeta)^{2(n+1+\alpha)+{\frac{1}{2}}}}dv_{\alpha}(z) + \int_{(2 Q)^{c}}\frac{r}{d(z,\zeta)^{2(n+1+\alpha)+1}}dv_{\alpha}(z).
\end{split}\ee
The second term on the right hand side has been estimated in the proof of Lemma \ref{le:AtEst}. The first term can be estimated as follows:
\be\begin{split}
\int_{d(z,\,\zeta)> 2 r} & \frac{r^{1/2} }{d(z,\,\zeta)^{2(n+1+\alpha)+1/2}} d v_{\alpha}(z) =r^{1/2} \sum_{k \ge 0} \int_{2^k r < d(z,\,\zeta) \le 2^{k+1} r } \frac{1}{d(z,\,\zeta)^{2(n+1+\alpha)+1/2}} d v_{\alpha}(z)\\
& \le r^{1/2} \sum_{k \ge 0} \frac{v_{\alpha} (Q_{2^{k+1} r})}{(2^k r)^{2(n+1+\alpha)+1/2}} \le C r^{1/2}\sum_{k=0}^\infty\frac{(2^{k+1} r)^{2(n+1+\alpha)}}{(2^k r)^{2(n+1+\alpha)+1/2}}\le C,
\end{split}\ee
where we have used the fact that $v_{\alpha} (Q_r) \approx r^{2(n+1+\alpha)}$ in the third inequality (cf. \cite[Corollary 5.24]{Zhu2005}).

By the same argument we can estimate $I_2$ and omit the details.

Next, we estimate $I_3.$ Note that
\be\begin{split}
M_3 \leq & (1-|\langle w,u\rangle|^2)|1-\langle w,\zeta\rangle|^{n+2+\alpha} \left | (1-\langle w,\zeta\rangle)^{n+2+\alpha}-(1-\langle w,u\rangle)^{n+2+\alpha} \right |\\
\leq & 2 |1-\langle w,u\rangle||1-\langle w,\zeta\rangle|^{n+2+\alpha} \left | \int^1_0\frac{d}{dt}(1-\langle w,t\zeta+(1-t)u\rangle)^{n+2+\alpha} d t \right |\\
= & 2 (n+2+\alpha) |1-\langle w,u\rangle||1-\langle w,\zeta\rangle|^{n+2+\alpha}\\
& \quad \times \left | \langle w,\zeta-u \rangle \int^1_0 (1-\langle w, t\zeta+(1-t) u \rangle )^{n+1+\alpha} d t \right |\\
\leq & C_{\g} |1-\langle w,u\rangle||1-\langle w,\zeta\rangle|^{n+2+\alpha}r |1-\langle
z,\zeta\rangle|^{n+3/2+\alpha},
\end{split}\ee
where the last inequality is achieved by the following estimates
\be\begin{split}
|1-\langle w,t\zeta+(1-t)u\rangle| & \leq C_{\g} |1-\langle z,t\zeta+(1-t)u\rangle|\\
& \leq C_{\g} |1-\langle z,u\rangle|+|\langle z,\zeta-u\rangle| \leq C_{\g} |1-\langle z,\zeta\rangle|
\end{split}\ee
and
\be
|\langle w, \zeta - u \rangle |\leq C_{\g} r|1-\langle z,\zeta\rangle|^{\frac{1}{2}},
\ee
for any $w \in D(z, \g)$ and $u \in Q_r (\zeta).$ Thus, by Lemmas \ref{le:KerEstimation} and \ref{le:KerEstimation2}
\be\begin{split}
I_{3} \leq & C_{\g} \left ( \int_{D(z,\gamma)} \sup_{u \in Q}\frac{(1-|w|^2)^{\frac{q}{2}} r^{\frac{q}{2}} |1-\langle z,\zeta\rangle|^{\frac{q}{2}(n+ \frac{3}{2} + \alpha)}}{|1-\langle w,u\rangle|^{q(n+1+\alpha) + \frac{q}{2}}|1-\langle w,\zeta \rangle|^{\frac{q}{2}(n+2+\alpha)}}d \tau (w) \right )^{\frac{1}{q}}\\
\leq & C_{\g} \left( \int_{D(z,\gamma)} \sup_{u \in Q} \frac{(1-|z|^{2})^{\frac{q}{2}} r^{\frac{q}{2}}|1-\langle z,\zeta\rangle|^{\frac{q}{2}(n+ \frac{3}{2}+\alpha)}}{|1-\langle
z,u\rangle|^{q(n+1+\alpha) + \frac{q}{2}}|1-\langle z,\zeta \rangle|^{\frac{q}{2}(n+2+\alpha)}} d \tau (w) \right )^{\frac{1}{q}}\\
\le & C_{\g} \Big ( \frac{(1-|z|^2)^{\frac{q}{2}} r^{\frac{q}{2}} }{| 1 - \langle z,\zeta \rangle|^{q (n+1+ \alpha)+ \frac{3}{4}q}} \Big )^{\frac{1}{q}}
\leq C_{\g} \frac{r^{\frac{1}{2}}}{d (z, \zeta)^{2(n+ 1+\alpha) + \frac{1}{2}}}.
\end{split}\ee
Hence,
\be
\int_{(2 Q)^{c}}I_3dv_{\alpha}(z)\leq C_{\g} \int_{(2 Q)^{c}} \frac{r^{\frac{1}{2}}}{d (z, \zeta)^{2(n+ 1+\alpha) + \frac{1}{2}}} d v_{\alpha}(z) \leq C_{\g},
\ee
as shown above.

Similarly, we can estimate $I_4$ and omit the details. Therefore, combining above estimates we conclude that
\be
\int_{(2 Q)^{c}} A_{\gamma} (\tilde{\nabla} P_{\alpha} a) d v_{\alpha} \leq C,
\ee
where $C$ depends only on $q, \g, n, $ and $\alpha.$
\hfill$\Box$

\begin{remark}\label{rk:atomp<1}\rm
We remark that whenever $\A^p_{\alpha}$ have an atomic decomposition in terms of atoms with respect to Carleson tubes for $0<p<1,$ the argument of Theorem \ref{th:AreaCharat} works as well in this case. However, as noted in Remark \ref{rk:atomdecomp-p<1} (1), the problem of the atomic decomposition of $\A^p_{\alpha}$ with respect to Carleson tubes in $0<p<1$ is entirely open.
\end{remark}

\begin{remark}\label{rk:AreaCharactp>1}\rm
The area integral inequality in case $1< p < \8$ can be also proved through using the method of vector-valued Calder\'{o}n-Zygmund operators for Bergman spaces. This has been done in \cite{CO2011}.
\end{remark}

\subsection*{Acknowledgement} This research was supported in part by the NSFC under Grant No. 11171338.

\end{document}